\documentstyle{amsppt}
\magnification 1200
\def\today
{\ifcase\month\or
     January\or February\or March\or April\or May\or June\or
     July\or August\or September\or October\or November\or December\fi
     \space\number\day, \number\year}
\magnification 1200
\input pictex
\UseAMSsymbols
\hsize 5.5 true in
\vsize 8.5 true in
\parskip=\medskipamount
\NoBlackBoxes

\def\mathbb{\Bbb}

\def\mathcal{\Cal}

\def\ve{\varepsilon}
\def\vp{\varphi}
\def\arrowk{^\to{\kern -6pt\topsmash k}}
\def\arrowK{^{^\to}{\kern -9pt\topsmash K}}
\def\arrowr{^\to{\kern-6pt\topsmash r}}
\def\bark{\bar{\kern-0pt\topsmash k}}
\def\arrowvp{^\to{\kern -8pt\topsmash\vp}}
\def\arrowf{^{^\to}{\kern -8pt f}}
\def\arrowg{^{^\to}{\kern -8pt g}}
\def\arrowu{^{^\to}a{\kern-8pt u}}
\def\arrowt{^{^\to}{\kern -6pt t}}
\def\arrowe{^{^\to}{\kern -6pt e}}
\def\tk{\tilde{\kern 1 pt\topsmash k}}
\def\barm{\bar{\kern-.2pt\bar m}}
\def\barN{\bar{\kern-1pt\bar N}}
\def\barA{\, \bar{\kern-3pt \bar A}}

\def\mathbb{\Bbb}

\def\snint{\raise2pt\hbox{$_{^\not}$}\kern-3.5 pt\int}
\def\nint{\not\!\!\int}
\def\dist{\text{\rm dist\,}}

\def\mes{\text {\rm mes\,}}
\TagsOnRight
\NoRunningHeads

\document
\topmatter
\title
Moment Inequalities for Trigonometric Polynomials with Spectrum in Curved Hypersurfaces
\endtitle
\author
J.~Bourgain
\endauthor
\address
Institute for Advanced Study, Princeton, NJ 08540
\endaddress
\endtopmatter

\noindent
{\bf (0). Summary}

In this note we develop further the technique from [B-G], based on the multi-linear restriction theory from [B-C-T], to establish some new
inequalities on the distribution of trigonometric polynomials on the $n$-dimensional torus $\Bbb T^n$, $n\geq 2$, of the form
$$
f(x) =\sum_{z\in\Cal E} a_z e^{2\pi ix.z}\tag 0.1
$$
where $\Cal E$ stands for the set of $\Bbb Z^n$-points on some dilate D.S of a fixed compact, smooth hypersurface $S$ in $\Bbb R^n$ with positive
definite second fundamental form.
More precisely, we prove that for $p\leq \frac {2n}{n-1}$ and any fixed $\ve>0$, the bound
$$
\Vert f\Vert_{L^p(\Bbb T^n)} \leq C_\ve D^\ve\Vert f\Vert_{L^2(\Bbb T^n)}\tag 0.2
$$
holds.

In particular, if $\Delta$ stands for the Laplacian on $\Bbb T^n$ and
$$
-\Delta f= Ef\tag 0.3
$$
we have that for $p\leq\frac {2n}{n-1}$, $n\geq 2$
$$
\Vert f\Vert_{L^p(\Bbb T^n)}\ll_\ve E^\ve \Vert f\Vert_{L^2(\Bbb T^n)}.\tag 0.4
$$
Recall that if $n=2$, one has the inequality, for $f$ satisfying (0.3),
$$
\Vert f\Vert_{L^4(\Bbb T^2)}\leq C\Vert f\Vert _{L^2(\Bbb T^2)}\tag 0.5
$$
due to Zygmund and Cook.
For $n=3$, arithmetical considerations permit to  obtain a bound
$$
\Vert f\Vert_{L^4(\Bbb T^3)}\ll_\ve E^\ve \Vert f\Vert_{L^2(\Bbb T^3)}\tag 0.6
$$
For $n\geq 4$, no estimate of the type (0.4) for some $p>2$ seemed to be known.
Recall also that it is {\it conjectured} that one has uniform bounds
$$
\Vert f\Vert_{L^q(\Bbb T^n)} \leq C_q \Vert f\Vert_{L^2(\Bbb T^n)} \text { if } q<\frac {2n}{n-2}\tag 0.7
$$
and
$$
\Vert f\Vert_{L^q(\Bbb T^n)} \leq C_q E^{\frac 12(\frac {n-2}2 -\frac nq)}\Vert f\Vert_{L^2(\Bbb T^n)} \text { if } q>\frac {2n}{n-2}
\tag 0.8
$$
if $f$ satisfies (0.3).
The inequality (0.8) was proven in [B1] (using the Hardy-Littlewood circle method) under the assumption
$$
q>\frac {2(n+1)}{n-3}\tag 0.9
$$
(up to an $E^\ve $-factor).

Another application of (0.2) relates to the periodic Schr\"odinger group $e^{it\Delta}$.
For $n\geq 1$, one has the Strichartz' type inequality
$$
\Vert (e^{it\Delta} f)(x)\Vert_{L^q(\Bbb T^{n+1})} \ll R^\ve \Vert f\Vert_{L^2(\Bbb T^n)}\tag 0.10
$$
for $q\leq\frac{2(n+1)} n$ and $f$ satisfying supp\,$\hat f\subset\mathbb Z^n\cap B(0, R)$.

Combined with results from [B3], (0.10) implies that for $q>\frac {2(n+3)}{n}$
$$
\Vert (e^{it\Delta} f)(x)\Vert_{L^q(\Bbb T^{n+1})} \leq C_q R^{\frac n2-\frac {n+2}q} \Vert f\Vert_{L^2(\Bbb T^n)}\tag 0.11
$$
for $f$ as above.  Note that inequality (0.11) is optimal.
This result is new (and of interest to the theory of the nonlinear Schr\"odinger equations with periodic boundary conditions) for $n\geq 4$. (See
[B3] for more details).

More generally, fix a smooth function $\psi:U\to\Bbb R$ on a neighborhood $U$ of $0\in\Bbb R^n$ such that $D^2\psi$ is positive definite.
For $q\leq\frac{2(n+1)} n$ and $R\to\infty$,
$$
\align
&\Big[\int_{[0, 1]^{n+1}}\Big|\sum_{z\in\Bbb Z^n, |z|<R} a_z e^{2\pi i(x.z+R^2 t\psi(\frac zR))}\Big|^q \ dxdt\Big]^{1/q}\\
&\ll R^\ve \Big(\sum|a_z|^2\Big)^{\frac 12}.\tag 0.12
\endalign
$$
Taking $\psi(x) =\alpha_1x_1^2+\cdots+\alpha_nx_n^2, \alpha_1, \ldots, \alpha_n>0$, generalizes (0.10) to irrational tori (cf. [B]).

\bigskip
\noindent
{\bf (1). Multilinear Estimates}

Fix a smooth, compact hyper-surface $S$ in $\Bbb R^n$ with positive definite second fundamental form.
For $x\in S$, denote $x' \in S^{(n-1)} = [|x|=1]$ the normal vector at the point $x$ and let $\sim :S^{(n-1)} \to S$ be the Gauss map.
Thus $\tilde {x'}=x$ for $x\in S$. 
Let $\sigma$ be the surface measure of $S$.

The estimates below depend on the multi-linear theory developed in [BCT] to bound oscillatory integral operators.
We recall the following version for later use.  Let
$$
\phi(x, y) = x_1y_1 +\cdots+ x_{n-1} y_{n-1} + x_n \big(\langle Ay, y\rangle +O(|y|^3)\big)\tag 1.1
$$
where $x\in\Bbb R^n$, $y\in\Bbb R^{n-1}$ is restricted to a small neighborhood of $0$ and $A$ is symmetric and definite (in particular, $A$ is non-degenerate).

Denote
$$
Z(x, y) =\partial_{y_1} (\nabla_x\phi)\wedge\cdots\wedge \partial_{y_{n-1}}(\nabla_x\phi).\tag 1.2
$$
Fix $2\leq k\leq n$ and disjoint balls $U_1,\ldots, U_k\subset\Bbb R^{n-1}$ such that the transversality condition holds
$$
|Z(x, y^{(1)})\wedge\cdots\wedge Z(x, y^{(k)})|> c\text { for all $x$ and $y^{(i)} \in U_i$}.\tag 1.3
$$
Then
$$
\Big\Vert\Big(\prod^k_{i=1} |T f_i|)^{\frac 1k}\Big\Vert_{L^q(B_R)} \ll R^\ve \Big(\prod^k_{i=1} \Vert f_i\Vert_2\Big)^{\frac 1k}\tag 1.4
$$
with $q=\frac{2k}{k-1}$, provided supp\,$f_i\subset U_i$.

\bigskip
\noindent
{\bf (2). Preliminary Lemmas}

We recall a few estimates from [B-G], \S3.

\proclaim
{Lemma 1}

Let $U_1, \ldots, U_n\subset S$ be small caps such that $|x_1'\wedge\cdots\wedge x_n'|>c$ for $x_i\in U_i$.

Let $M$ be large and $\Cal D_i\subset U_i (1\leq i\leq n)$ discrete sets of $\frac 1M$-separated points.

Let $B_M\subset \Bbb R^n$ be a ball of radius $M$. Then, for $q=\frac {2n}{n-1}$
$$
\nint_{B_M} \prod^n_{i=1} \Big|\sum_{\xi\in\Cal D_i} a(\xi) e^{ix.\xi}\Big|^{q/n}\ll M^\ve \prod^n_{i=1} \Big[\sum_{\xi\in \Cal D_i}
 |a(\xi)|^2\Big]^{\frac q{2n}}\tag 2.1
$$
where $\snint$ denotes the average.
\endproclaim

\noindent
{\bf Proof.}

This is just a discretized version of (2.4) with $k=n$; our assumption ensures the required transversality condition (1.3)

We can assume $B_M$ centered at $0$.
Introduce functions $g_i$ on $U_i$ defined by
$$
\Bigg\{
\aligned g_i(\zeta)&= a(\xi) \text { if } |\zeta-\xi|< \frac cM, \xi\in\mathcal D_i\\
g_i(\zeta)&= 0 \text { otherwise}.\endaligned
 \tag 2.2
$$
($c>0$ a small constant). One may then replace $\sum_{\xi\in\mathcal D_i} a(\xi) e^{ix. \xi}$ by
\break $c'M^{n-1} \int_S g_i(\zeta) e^{ix.\zeta} \sigma (d\zeta)\text { if } x\in B_M$.
Hence
$$
\align
&\int_{B_M} \prod^n_{i=1} \Bigg|\sum_{\zeta \in\mathcal D_i} a(\xi) e^{ix.\xi}\Bigg|^{q/n} dx \lesssim\\
&M^{(n-1)q} \int_{B_M} \prod^n_{i=1} \Big|\int_S g_i(\zeta) e^{ix\zeta} \sigma(d\zeta)\Big|^{q/n} dx \overset{\text
{(1.4)}}\to \ll\\
&M^{(n-1)q+\ve}\prod^n_{i=1} \Vert g_i\Vert^{q/n}_{L^2(U_i)} \sim M^{\frac {n-1}2 q+\ve} \prod^n_{i=1}
\Bigg[\sum_{\xi\in\mathcal D_i}|a(\xi)|^2\Bigg]^{\frac
q{2n}}.\tag 2.3  \endalign
$$
Since $\snint_{B_M}$ refers to the average, (2.1) follows, since $q=\frac {2n}{n-1}$.

\proclaim
{Lemma 2}

 Let $S\subset\mathbb R^n$ be as above and $2\leq m\leq n$.
Let $V$ be an $m$-dimensional subspace of $\mathbb R^n$, $P_1, \ldots, P_m\in S$ such that
$$
P_1', \ldots, P_m' \in V \text { and } |P_1\wedge\cdots\wedge P_m|>c \tag 2.4
$$
and $U_1, \ldots, U_m \subset S$ sufficiently small neighborhoods of $P_1, \ldots, P_m$.

Let $M$ be large and $\mathcal D_i\subset U_i$ $(1\leq i\leq m)$ discrete sets of $\frac 1M$-separated points $\xi\in
S$ such that $\dist(\xi', V)<\frac cM$.
Let $g_i\in L^\infty (U_i)(1\leq i\leq m)$.
Then letting $q=\frac {2m}{m-1}$
$$
\aligned
&\nint_{B_M} \prod^m_{i=1} \Bigg|\sum_{\xi\in\mathcal D_i}\Big(\int_{|\zeta -\xi|<\frac cM} g_i(\zeta) e^{ix.\zeta}
\sigma(d\zeta)\Big)\Bigg|^{q/m} dx\ll\\
&M^\ve\Big\{\nint_{B_M} \prod^m_{i=1} \Bigg[\sum_{\xi\in\mathcal D_i} \Big|\int_{|\zeta-\xi|<\frac cM} g_i(\zeta)
e^{ix.\zeta} \sigma(d\zeta)\Big|^2\Bigg]^{1/2m}\Bigg\}^q.
\endaligned
\tag 2.5
$$
\bigskip
\endproclaim

\noindent
{\bf Proof.}

Performing a rotation, we may assume $V=[e_1, \ldots, e_m]$ and denote $\tilde V\subset S$ the image of
$V\cap S^{(n-1)}$ under the Gauss map.
Let again $B_M$ be centered at $0$.
For each $\xi\in \bigcup^m_{i=1}\mathcal D_i$ there is by assumption some $\hat\xi\in \tilde V$.
$|\xi-\hat\xi|< \frac cM$.
Write
$$
\int_{|\zeta-\xi|<\frac cM} g_i(\zeta) e^{ix.\zeta} \sigma(d\zeta)=e^{ix\hat \xi}\int_{|\zeta-\xi|<\frac cM}
g_i(\zeta) e^{ix.(\zeta-\hat\xi)} \sigma (d\zeta).\tag 2.6
$$
Since in the second factor of (2.6), $|\zeta-\hat\xi|=o(\frac 1M)$, we may view it as constant $a(\xi)$ on $B_M\subset \mathbb
R^n$.

Thus we need to estimate
$$
\nint_{B_M} \Big\{\prod^m_{i=1} \Big|\sum_{\xi\in\Cal D_i} e^{ix.\hat\xi}a(\xi)\Big|^{q/m}\Big\} dx.\tag 2.7
$$
Writing $x=(u, v)\in B_M^{(m)} \times B_M^{(n-m)}$, (2.7) may be bounded by
$$
\max_{v\in B_M^{(n-m)}} \nint_{B_M^{(m)}} \Big\{\prod^m_{i=1} \Big\}\Big|\sum_{\xi\in\Cal D_i} e^{iu.\pi_m(\hat\xi)}
a_v(\xi)\big|^{q/m}\Big\} du\tag 2.8
$$
with $a_v(\xi)= e^{iv.\hat\xi} a(\xi)$.

Since $S$ has positive definite second fundamental form, $\pi_m(\tilde V)\subset V=[e_1, \ldots, e_m]$ is a
hypersurface in $V$ with same property and the normal vector at $\pi_m(\hat\xi)=(\hat\xi)' \in V$.
Since (2.4), application of (2.1) with $n$ replaced by $m$ and $\Cal D_i$ by $\{\pi_m\hat\xi; \xi \in\Cal D_i\}$ gives
the estimate on (2.7)
$$
\ll M^\ve\prod^m_{i=1} \Big[\sum_{\xi\in\Cal D_i} |a(\xi)|^2\Big]^{q/2m}
$$
and (2.5) follows .

\proclaim
{Lemma 3}
Let
$$
p=\frac {2n}{n-1}.
$$
Take $K_n\gg K_{n-1}\gg \cdots\gg K_1\gg 1$.
For $1\leq j\leq n$, denote by $\{U_\alpha^{(j)}\}$ a partition of $S$ in cells of size $\frac 1{K_j}$.
Then, for $R>K_n$ and $g\in L^2(S)$,
$$
\spreadlines{6pt}
\align
\Big\Vert\int g(\xi) e^{ix.\xi} \sigma(d\xi)\Big\Vert_{L^p(B_R)} &\ll_\ve\\
C(K_n) R^\ve\Big[\int_S|g(\xi)|^2 \sigma(d\xi)\Big]^{\frac 12}&+\sum_{2\leq j\leq n} C(K_{j-1})K_j^\ve
\Big\{\sum_\alpha \Big\Vert\int_{U_\alpha(j)} g(\xi) e^{ix.\xi}\sigma (d\xi)\Big\Vert^2_{L^p(B_R)}\Big\}^{\frac 12}\\
&+\Big\{\sum_\alpha\Big\Vert\int_{U_\alpha^{(1)}} g(\xi) e^{ix\xi} \sigma(d\xi)
\Big\Vert^2_{L^p(B_R)}\Big\}^{\frac 12}\tag 2.9
\endalign
$$
where $C(K)$ denotes some polynomial function of $K$.
\endproclaim

\noindent
{\bf Proof.}
We follow the analysis from \S3 in [B-G].

For $x\in B_R$, let
$$
(2.10) =\int_Sg(\xi) e^{ix.\xi}\sigma (d\xi)
$$

Start decomposing $S=\bigcup_\alpha U_\alpha (\frac 1{K_n})$ in caps of size $\frac 1{K_n}$ and write
$$
(2.10) =\sum_\alpha\int_{U_\alpha(\frac 1{K_n})} g(\xi) e^{ix.\xi}\sigma (d\xi)=\sum_\alpha c_\alpha(x).
$$
Fixing $x$, there are 2 possibilities

(2.11) There are $\alpha_1, \alpha_2, \ldots, \alpha_n$ such that
$$
|c_{\alpha_1}(x)|, \ldots, |c_{\alpha_n}(x)|> K_n^{-(n-1)}\max_\alpha |c_\alpha(x)|\tag 2.12
$$
and
$$
|\xi_1\wedge\cdots\wedge \xi_n| \gtrsim K_n^{-n} \text { for } \xi_i\in U_{\alpha_i}.
\tag 2.13
$$
(2.14)
\  The negation of (2.11), which implies that there is an $(n-1)$-dim subspace $V_{n-1}$ such that
$$
|c_\alpha(x)|\leq K_n^{-(n-1)} \max_\alpha |c_\alpha(x)|\ \text { if } \ \dist (U_\alpha, \tilde V_{n-1})\gtrsim \frac 1{K_n}.
$$
If (2.11), it follows from (2.12) that
$$
\Big|\int_Sg(\xi) e^{ix.\xi} \sigma (d\xi)\Big| \leq K_n^{n-1} \max |c_\alpha(x)| \leq K_n^{2n-2} \Big[\prod^n_{i=1}
|c_{\alpha_i}(x)|\Big]^{\frac 1n}
$$
and the corresponding contribution to the $L^p_{B_R}$-norm of (4.1) is bounded by
$$
\align
\int^{(2.11)}_{B_R}&  \Big|\int_Sg(\xi) e^{ix.\xi} \sigma(d\xi)\Big|^p\\
& \lesssim K_n^{2p(n-1)} \sum_{\Sb \alpha_1,
\ldots, \alpha_n\\ {(2.13)}\endSb}
\int_{B_R}\prod^n_{i=1} \Big|\int_{U_{\alpha_i(\frac 1{K_n})}} g(\xi) e^{ix.\xi}\sigma(d\xi)\Big|^{\frac pn}.
\tag 2.15
\endalign
$$
In view of (2.13), the \cite{BCT}-estimate (1.4) with $k=n$ applies to each (2.15) term.
Thus
$$
\int_{B_R} \prod^n_{i=1} \Big|\int_{U_{\alpha_i(\frac 1{K_n})}} g(\xi) e^{ix.\xi} \sigma (d\xi)\Big|^{\frac 2{n-1}}
dx\ll C(K_n)R^\ve \Big[\int_S|g(\xi)|^2 \sigma(d\xi)\Big]^{\frac n{n-1}}.\tag 2.16
$$
Next consider the case (2.14).Thus
$$
\aligned
&|(2.10)|\leq\Big|\int_{\dist(\xi, \tilde V_{n-1})\lesssim \frac 1{K_n}} g(\xi) e^{ix.\xi}\sigma(d\xi)\Big|+
\max_\alpha \Big|\int_{ U_\alpha (\frac 1{K_n})} g(\xi) e^{ix.\xi}\sigma (d\xi)\Big|\\
&=(2.17)+(2.18)
\endaligned
$$
where $V_{n-1}$ depends on $x$.

Note however that, from its definition, we may view $|c_\alpha(x)|$ as `essentially' constant on balls of size $K_n$.
Making this claim rigorous requires some extra work and one replaces $|c_\alpha(x)|$ by a majorant $|c_\alpha|*
\eta_{K_n}$, $\eta_K(x)=\frac 1{K^n} \eta\Big(\frac xK\Big)$ and $\eta$ a suitable bump-function.
We may then ensure that $|c_\alpha|*\eta _{K_n}$ is approximately constant at scale $K_n$.
But we will not sidetrack the reader with these technicalities that may be found in [B-G], \S2.

Thus, upon viewing the $|c_\alpha|$ approximatively constant at scale $K_n$, the bound (2.17) + (2.18) may clearly be
considered valid on $B(\bar x, K_n)$ with the same linear space $V_{n-1}$.

Obviously
$$
(2.18) \leq \Big(\sum_\alpha \Big|\int_{U_\alpha (\frac 1{k_n})} g(\xi) e^{ix.\xi}\sigma(d\xi)|^p\Big)
^{\frac 1p}
$$
and the corresponding $L^p_{B_R}$-contribution is bounded by
$$
\Big\{\sum_\alpha \Big\Vert\int_{U_\alpha (\frac 1{K_n})} g(\xi) e^{ix.\xi} \sigma(d\xi) \Big\Vert^2_{L^p_{B_R}}
\Big\}^{1/2}.\tag 2.19
$$
Consider the term (2.17).
Proceeding similarly, write for $x\in B(\bar x, K_n)$
$$
\aligned
\int_{\dist (\xi, V_{n-1})\lesssim \frac 1{K_n}} g(\xi) e^{ix.\xi} \sigma (d\xi)&=\\
\sum_\alpha\int_{U _\alpha (\frac 1{K_{n-1}})\cap[\dist (\xi, \tilde V_{n-1})\lesssim \frac 1{K_n}]} &g(\xi)
e^{ix.\xi}\sigma(d\xi) =\sum_\alpha
c_\alpha^{(n-1)}(x).
\endaligned
\tag 2.20
$$
We distinguish the cases

(2.20) There are $\alpha_1, \ldots, \alpha_{n-1}$ such that
$$
|c_{\alpha_1}^{(n-1)} (x)|, \ldots, |c^{(n-1)}_{\alpha_{n-1}} (x)|>
K_{n-1}^{-(n-2)}\max_\alpha|c_\alpha^{(n-1)}(x)|\tag 2.21
$$
and
$$
|\xi_1'\wedge\ldots\wedge\xi_{n-1}'|\gtrsim  K_{n-1}^{-(n-1)}  \ \text { for } \ \xi_i\in U_{\alpha_i} \Big(\frac
1{K_{n-1}}\Big).\tag 2.22
$$

(2.23) Negation of (2.20), implying that there is an $(n-2)$-dim subspace $V_{n-2}\subset V_{n-1}$ (depending on $x$)
such that
$$
|c_\alpha^{(n-1)}(x)|< K_{n-1}^{-(n-2)} \max_\alpha |c_\alpha^{(n-1)} (x)| \text { for } \dist (U_\alpha,
\tilde V_{n-2})\gtrsim \frac 1{K_{n-1}}.
$$
This space $V_{n-2}$ can then again be taken the same on a $ K_{n-1}$-neighborhood of $x$.

We analyze the contribution of (2.20).  By (2.21)
$$
|(2.19)| < K_{n-1}^{2n-4} \Big[\prod^{n-1}_{i=1} |c_{\alpha_i}^{(n-1)} (x) |\Big]^{\frac 1{n-1}}\tag 2.24
$$
and hence
$$
\aligned
&\operatornamewithlimits\nint\limits_{\Sb B(\bar x, K_n)\\ x\text { satisfies }
(2.20)\endSb}\Big|\int\limits_{\dist(\xi, \tilde V_{n-1})\lesssim \frac 1{K_n}} g(\xi) e^{ix.\xi}
\sigma(d\xi)\Big|^p\leq\\
&K_{n-1}^{p(2n-4)}\sum_{\Sb {\alpha_1, \ldots, \alpha_{n-1}}\\ {(2.22)}\endSb} \nint_{B(\bar x,
K_n)}\Big\{\prod^{n-1}_{i=1} \Big|\int_{U_{\alpha_i}
(\frac 1{K_{n-1}}) \cap [\dist(\xi, \tilde V_{n-1})\lesssim \frac 1{K_n}]} g(\xi) e^{ix.\xi} \sigma(d\xi)\Big|^{p/n-1}\Big\}.
\endaligned
\tag 2.25
$$
We use the bound (2.5) to estimate the individual integrals
$$
(2.26) \ \nint_{B(\bar x, K_n)} \Big\{\prod^{n-1}_{i=1} \Big|\int_{U_{\alpha_i} (\frac 1{K_{n-1}}) \cap [\dist(\xi,
\tilde V_{n-1})\lesssim \frac 1{K_n}]}
g(\xi)e^{ix.\xi} \sigma(d\xi) \Big|\Big\}^{\frac q{n-1}} \text { with } q=\frac{2(n-1)}{n-2}.
$$
Thus $m=n-1, V=V_{n-1}$ and $P_i$ is the center of $U_{\alpha_i}(\frac 1{K_{n-1}})$.
Let $M=K_n$ and $\Cal D_i$ the centers of a cover of $U_{\alpha_i} (\frac 1{K_{n-1})} \cap [\dist (\xi, \tilde V_{n-1})\lesssim \frac 1{K_N}]$ by caps $U_\alpha(\frac
1{K_n})$.

By (2.5) we get an estimate 
$$
(2.26) \ll K_n^\ve C(K_{n-1}) \Big\{\nint_{B(\bar x, K_n)} \prod^{n-1}_{i=1} \Big[\sum_\alpha^{(i)}
\Big|\int_{U_\alpha(\frac 1{K_n})}
g(\xi) e^{ix.\xi} \sigma(d\xi)\Big|^2\Big]^{\frac 1{2(n-1)}}\Big\}^q\eqno{(2.27)}
$$
where in $\sum^{(i)}$ the sum is over those $\alpha$ such that $U_\alpha(\frac 1{K_n})\subset U_{\alpha_i}(\frac
1{K_{n-1}})$ and
$U_\alpha(\frac 1{K_n})\cap \tilde V_{n-1} \not=\phi$.
Hence, we certainly have
$$
(2.26)\ll K_n^\ve \ C(K_{n-1})\Big\{\nint_{B(\bar x, K_n)} \Bigg[\sum_{\alpha}\Big|\int_{U_\alpha(\frac 1{K_n})}
g(\xi) e^{ix.\xi}\sigma(d\xi)\Big|^2 \Bigg]^{\frac 12}\Big\}^q
$$
and therefore, since $p<q$,
$$
(2.25) \ll K^\ve_nC(K_{n-1}) \Big\{\nint_{B(\bar x, K_n)} \Big[\sum_\alpha \Big|\int_{U_\alpha(\frac 1{K_n})}
g(\xi) e^{ix.\xi} \sigma(d\xi)\Big|^2\Big]^{p/2}\Big\}.\tag 2.28
$$
Hence the collected contribution over $B_R$ of (2.28) is bounded by
$$
K_n^\ve C(K_{n-1})\Big\{\sum_\alpha \Big\Vert\int_{U_\alpha(\frac 1{K_n})} g(\xi) e^{ix.\xi} \sigma(d\xi)
\Big\Vert^2_{L^p(B_R)}\Big\}^{1/2}.\tag 2.29
$$

Next, we analyze the contribution of (2.23) which is similar to that of (2.14) with $n-1$ replaced by $n-2$ and
$K_n$ by $K_{n-1}$.
The local estimate (2.27) becomes
$$
K^\ve_{n-1} C(K_{n-2})\Big\{\nint_{B(\bar x, K_{n-1})} \prod^{n-2}_{i=1} \Big[\sum^{(i)}_{\alpha}
\big|\int_{U_\alpha(\frac 1{K_{n-1}})} g(\xi) e^{ix.\xi}\sigma(d\xi)\Big|^2\Big]^{\frac 1{2(n-2)}}\Big\}^q
\tag 3.30
$$
with $q= \frac{2(n-2)}{n-3}$ and where in $\sum^{(i)}$ the sum is over those $\alpha$ such that
$$
U_\alpha \Big(\frac 1{K_{n-1}}\Big) \subset U_{\alpha_i}\Big(\frac 1{K_{n-2}}\Big) \text { and }
U_\alpha\Big(\frac 1{K_{n-1}}\Big) \cap \tilde V_{n-2} \not= \phi.
$$

The collected contribution of (2.30) to the $L^p_{B_R}$-norm of (2.10) is bounded by
$$
K^\ve_{n-1} C(K_{n-2}) \Big\{\sum_\alpha \Big\Vert \int_{U_\alpha(\frac 1{K_{n-1}})} g(\xi) e^{ix.\xi}
\sigma(d\xi)\Big\Vert^2_{L^p_{(B_R)}}\Big\}^{\frac 12}.\tag 3.31
$$
The continuation of the process is now clear and leads to the bound (2.9).  This proves Lemma 3.

\bigskip

Taking $K_j>K_{j-1}^{C/\ve}$ in Lemma 3, we obtain

\proclaim
{Lemma 4}
Fix $\ve>0$.
Let $K_1\gg 1$ be large enough and assume $R>K_1^{C(\ve)}$.

Then, with $p=\frac{2n}{n-1}$
$$
\align
\Big\Vert \int g(\xi) &e^{ix.\xi} \sigma(d\xi)\Big\Vert_{L^p_{B_R}}\leq R^\ve\Big[\int_S |g(\xi)|^2 \sigma
(d\xi)\Big]^{\frac 12}\\
&+\max_{K_1<K<K_1^{C(\ve)}} \Big\{K^\ve\sum_\alpha\Big\Vert\int_{U_\alpha (\frac 1K)}
g(\xi)e^{ix.\xi} \sigma(d\xi)\Big\Vert ^2_{L^p_{B_R}}\Big\}^{1/2}\tag 2.32
\endalign
$$
with $\{U_\alpha (\frac 1K)\}$ a cover of $S$ by $\frac 1K$-size caps.
\endproclaim

The first term on the right side of (2.32) may be eliminated.

Observe first that since $|x|<R$, the left side may be replaced by
$$
\Big\Vert\int G(\xi) e^{ix.\xi} \sigma(d\xi) \Big\Vert_{L^p_{B_R}}\tag 2.33
$$
where $G$ is a smoothing of $g$ at scale $\frac 1R$.

Applying (2.32) with $g$ replaced by $G$, the first term on the right
$$
\Big[\int_S|G(\xi)|^2 \sigma(d\xi)\Big]^{\frac 12} \lesssim \Big\{\sum_\alpha \Big\Vert\int_{U_\alpha(\frac cR)}
g(\xi) e^{ix.\xi} \sigma(d\xi)\Big\Vert^2_{L^p_{B_R}} \Big\}^{\frac 12}\tag 2.34
$$
and the other terms may be majorized by
$$
\Big\Vert\int_{U_\alpha(\frac 1K)} G(\xi) e^{ix.\xi}\sigma(d\xi)\Big\Vert_{L^p_{B_R}}\lesssim
\Big\Vert \int_{U_\alpha(\frac 1K)} g_1(\xi) e^{ix.\xi}\sigma(d\xi)\Big\Vert_{L^p_{B_R}}\tag 2.35
$$
for some $g_1=\eta g$ with $\eta$ a smooth function.

Hence we obtain

\proclaim
{Lemma 5}
Fix $\ve>0$.
Let $K_1\gg 1$ be large enough and assume $R>K_1^{C(\ve)}$.
Then, with $p=\frac {2n}{n-1}$, we have
$$
\align
\Big\Vert\int g(\xi) e^{ix.\xi}\sigma(d\xi)\Big\Vert_{L^p_{B_R}} &< R^\ve
\Big\{\sum_\alpha\Big\Vert\int_{U_\alpha(\frac cR)} g(\xi) e^{ix.\xi} \sigma(d\xi) \Big\Vert^2_{L^p_{(R)}}\Big\}^{\frac
12}+\tag 2.36\\
& \max_{K_1< K<K_1^{C(\ve)}} \Big\{K^\ve \sum_\alpha\Big\Vert\int_{U_\alpha(\frac 1K)} g(\xi) e^{ix.\xi}
\sigma(d\xi) \Big\Vert^2_{L^p_{(R)}}\Big\}^{\frac 12}
\endalign
$$
where $L^p_{(R)} =L^p\big(\omega(\frac xR) dx\big)$ with $0<\omega<1$ some rapidly decaying function on $\Bbb R^n$.
\endproclaim

\bigskip

In order to iterate (2.36), we rely on rescaling.

Parametrize $S$ (locally, after affine coordinate change) as
$$
\spreadmatrixlines{8pt}
\left\{
\matrix
\xi_i= y_i (1\leq i\leq n-1)\quad\qquad\\
\xi_n= y_1^2+\cdots+ y^2_{n-1} +O(|y|^3)
\endmatrix
\right.
\tag 2.37
$$
with $y$ taken in a small neighborhood of $0$.

Let $U(\rho)$ be a $\rho$-cap on $S$ and evaluate
$$
\Big\Vert\int_{U(\rho)} g(\xi) e^{ix.\xi} \sigma(d\xi) \Big\Vert_{L^p(B_R)}.\tag 2.38
$$
Thus in view of (2.37), (2.38) amounts to
$$
\Big\Vert\int_{B(a, \rho)} g(y) e^{i\vp(x, y)} dy\Big\Vert_{L^p(B_R)}\tag 2.39
$$
with
$$
\vp(x, y) =x_1y_1+\cdots+ x_{n-1} y_{n-1} +x_n \big(|y|^2+O(|y|^3)\big)
\tag 2.40
$$
and $B(a, \rho)\subset\Bbb R^{n-1}$.

A shift $y\mapsto y-a$ and change of variables $x_i'=x_i+x_n(2a_i+\cdots) \ (1\leq i<n)$ permits to set $a=0$.
By parabolic rescaling
$$
y=\rho y' \text { and } \rho x_i = x_i' (1\leq i<n), \rho^2 x_n=x_n'\tag 2.41
$$
we obtain a new phase function $\psi(x', y')$ and (2.39) becomes
$$
\rho^{n-1-\frac {n+1}p} \Big\Vert \int_{B(0, 1)} g(a+\rho y') e^{i\psi(x', y')} dy'\Big\Vert _{L^p(\Omega)}\tag 2.42
$$
where $\Omega =[|x_i'|< \rho R(1\leq i<n) , |x_n'| <\rho^2R]$.

Partition $\Omega =\bigcup \Omega_s$ in size-$\rho^2R$ balls $\Omega_s$ and apply Lemma 5 on each $\Omega_s$
with $R$ replaced by $\rho^2R$.
Assuming
$$
R>\rho^{-2} K_1^{C(\ve)}\tag 2.43
$$
(2.36) implies that
$$
\spreadlines{6pt}
\align
&\Big\Vert\int_{B(0,1)} g(a+\rho y') e^{i\psi(x', y')} dy'\Big\Vert_{L^p(\Omega_s)}<\\
&(\rho^2 R)^\ve \Big\{\sum_\alpha\Big\Vert\int_{U_\alpha (\frac c{\rho^2R})} g(a+\rho y') e^{i\psi(x', y')}
dy'\Big\Vert^2_{L^p(\omega(\frac {x'- b_s}{\rho^2R})dx')}\Big\}^{\frac 12}+\\
&\max_{K_1<K<K_1^{C(\ve)}} K^\ve\Big\{\sum_\alpha \Big\Vert\int_{U_\alpha(\frac 1K)} g(a+\rho y') e^{i\psi(x', y') }
dy'\Big\Vert^2_{L^p(\omega(\frac {x'-b_s}{\rho^2R}) dx')}\Big\}^{\frac 12}\qquad\tag 2.44
\endalign
$$
with $b_s$ the center of $\Omega_s$.

Note that certainly
$$
\sum_s\omega\Big(\frac {x'-b_s}{\rho^2 R}\Big)<\omega_1 \Big(\frac xR\Big).
$$
Summing $(2.44)^p$ over $s$ and reversing the coordinate changes clearly implies that
$$
\spreadlines{6pt}
\align
&(2.39), (2.42)<\\
&(\rho^2R)^\ve\Big\{\sum_\alpha \Big\Vert\int_{U_\alpha(\frac c{\rho R})} g(y) e^{i\vp (x, y)} dy\Big\Vert^2_{L^p_{(R)}}
\Big\}^{\frac 12}+\\
&\max_{K_1<K<K_1^{C(\ve)}} \Big\{ K^\ve\sum_\alpha\Big\Vert\int_{U_\alpha(\frac\rho K)} g(y) e^{i\vp(x, y)}
dy\Big\Vert^2_{L^p_{(R)}}\Big\}^{\frac 12}\tag 2.45
\endalign
$$
under the assumption (2.43).

Taking $R=\rho^{-2} K_2$ with $K_2> K_1^{C(\ve)}$ in (2.45), we obtain

\proclaim
{Lemma 6}
Let $K_2>K_1^{C(\ve)}$. Then
$$
\align
\Big\Vert&\int_{U(\rho)} g(\xi) e^{ix.\xi} \sigma(d\xi)\Big\Vert_{L^p(B_{K_2\rho^{-2})}}\\
&\ll_\ve \max_{K_1<K<K_2}\Big\{ K^\ve\sum_\alpha\Big\Vert\int_{U_\alpha(\frac {c\rho}K)}
g(\xi) e^{ix.\xi} \sigma(d\xi)\Big\Vert^2_{(K_2\rho^{-2})}\Big\}^{\frac 12}.\tag 2.46
\endalign
$$
\endproclaim

If $R>K_2\rho^{-2}$, we can partition $B_R$ in cubes of size $K_2\rho^{-2}$ and apply (2.46) on each of them,
with $g(\xi)$ replaced by $g(\xi)$ $e^{ia.\xi}$ for some $a\in B_R$.
Hence

\proclaim
{Lemma 6$'$}
Let $R>K_2\rho^{-2}, K_2=K_1^{C(\ve)}$. Then
$$
\align
\Big\Vert&\int_{U(\rho)} g(\xi) e^{ix.\xi} \sigma(d\xi)\Big\Vert_{L^p(B_R)}\\
&\ll_\ve \max_{K_1<K<K_2}\Big\{ K^\ve\sum_\alpha\Big\Vert\int_{U_\alpha(\frac {c\rho}K)}
g(\xi) e^{ix.\xi)} \sigma(d\xi)\Big\Vert^2_{L^p_{(R)}}\Big\}^{\frac 12}.\tag 2.47
\endalign
$$
\endproclaim

It is now straightforward to iterate Lemma $6'$ and derive the following statement

\proclaim
{Proposition 1} Let $0<\delta\ll 1$ and $R>C(\ve)\delta^{-2}$.
Then, with $p=\frac {2n}{n-1}$
$$
\Big\Vert\int g(\xi) e^{ix.\xi} \sigma(d\xi)\Big\Vert_{L^p_{(R)}} \ll_\ve
\delta^{-\ve}\Big\{\sum_\alpha \Big\Vert \int_{U_\alpha(\delta)} g(\xi) e^{ix.\xi}\sigma(d\xi)\Big\Vert^2_{L^p_{(R)}}
\Big\}^{\frac 12}.
\tag 2.48
$$
\endproclaim

\noindent
\hbox{\bf (3). $L^p$-bounds for certain exponential polynomials and applications}

We fix a smooth compact hyper-surface $S$ in $\Bbb R^n$ with positive definite second fundamental form.
We consider exponential polynomials with frequencies on some dilate $D.S$ of $S$.

\proclaim
{Proposition 2}
Let $0<\rho<D$ and let $\Cal E$ be a discrete set of points on the dilate $D.S$ that are mutually at least $\rho$
separated.
Then, for $p=\frac {2n}{n-1}$ and any (fixed) $\ve>0$
$$
\Big[\nint_{B_R} \Big|\sum_{z\in \Cal E} a_z e^{ix.z)} \Big\vert^p dx\Big]^{\frac 1p} \ll_\ve \Big(\frac
D\rho\Big)^\ve \Big(\sum_{z\in\Cal E}|a_z|^2\Big)^{\frac 12}\tag 3.1
$$
provided
$$
R>C(\ve)D{\rho^{-2}}.\tag 3.2
$$
\endproclaim

\noindent
{\bf Proof.}

By rescaling, we may clearly assume $D=1$.

Let $0<\tau <\rho /10$ and let $g$ be the function on $S$ defined by
$$
\aligned
g(\xi) &= \frac {a_z}{\sigma(U(z, \tau))} \text { if } \xi\in U(z, \tau)\\
&= 0 \ \text { otherwise}
\endaligned
\tag 3.3
$$
Here $U(z,\tau) \subset S$ denotes a $\tau$-neighborhood of $z$ on $S$.
Thus
$$
\int g(\xi) e^{ix.\xi} \sigma(d\xi) =\sum_{z\in\Cal E} a_z\nint_{U(z, \tau)} e^{ix.\xi} \sigma(d\xi).\tag 3.4
$$

Applying (2.48) with $\delta=\rho$, it follows from (3.3), (3.4) that
$$
\Big\{\nint_{B_R}\Big|\sum_{z\in\Cal E} a_z \nint_{U(z, \tau)} e^{ix.\xi} \sigma(d\xi)\Big|^p dx\Big\}
^{\frac 1p} \ll_\ve \rho^{-\ve}\Big(\sum_z |a_z|^2\Big)^{1/2}\tag 3.5
$$
letting $\tau\to 0$, (3.1) clearly follows.
\bigskip

Next, observe that if $\Cal E$ is contained in a lattice, then $\sum_{z\in\Cal E} a_z e^{ix.\xi}$ is a periodic
function.
Hence Proposition 2 implies

\proclaim
{Proposition 3}
Let $S$ be as above and $\Cal E=\Bbb Z^n \cap DS$, $D\to\infty$.

Then, with $p=\frac {2n}{n-1}$
$$
\Big[\int_{\Bbb T^n} \Big|\sum_{z\in\Cal E} a_z e^{2\pi ix.z}\Big|^p dx\Big]^{\frac 1p} \ll_\ve
D^\ve\Big(\sum|a_z|^2\Big)^{1/2}\tag 3.6
$$
where $\Bbb T^n$ stands for the $n$-dimensional torus.
\endproclaim

\proclaim
{Corollary 4} Let $\vp=\vp_E, -\Delta\vp_E=E\vp_E$ be an eigenfunction of $\Bbb T^n$, $n\geq 2$.
Then for $p=\frac {2n}{n-1}$ and any $\ve>0$, we have
$$
\Vert\vp\Vert_{L^p (\Bbb T^n)} \leq C(\ve)E^\ve \Vert\vp\Vert_{L^2(\Bbb T^n)}.\tag 3.7
$$
\endproclaim

\noindent
{\bf Remark.}
Corollary 4 should be compared with the result from [B1].
It is conjectured that for eigenfunctions of $\Bbb T^n, n\geq 2$, there is a uniform bound
$$
\Vert\vp\Vert_p\leq C(p) \Vert\vp\Vert_2\text { for } p<\frac {2n}{n-2}.\tag 3.8
$$
If $n=2$, (3.8) is known to hold for $p\leq 4$ (due to Zygmund-Cook) but for no exponent $p>4$.

If $n=3$, (3.7) is valid for $p\leq 4$.
This is a consequence of the following observation.
One clearly has the estimate
$$
\Vert\vp\Vert_4 \leq K^{1/4} \Vert\vp\Vert_2
$$
denoting
$$
K=\max_{\xi\in\ Z^3}\big(\#\{(\xi_1, \xi_2)\in\Bbb Z^3\times\Bbb Z^3; |\xi_1|^2 = E=|\xi_2|^2  \text { and
}\xi_1+\xi_2=\xi\}\big).
$$
Projecting on one of the coordinate planes reduces the  issue to bounding the number $|\Cal E\cap \Bbb Z^2|$ with $\Cal E\subset\Bbb R^2$
some ellipse of size at most $E^{1/2}$.
It is well known that
$$
|\Cal E\cap \Bbb Z^2|\ll E^{\ve}\tag 3.9
$$
(cf. [B-R]) and hence $K\ll E^\ve$.

For $n\geq 4$, no estimates of the type (3.7) for some $p>2$, seemed to be previously known.
Recall that for $n\geq 4$ and $R$ a large positive integer
$$
|RS^{(n-1)} \cap\Bbb Z^n|\sim R^{n-2}.\tag 3.10
$$
Thus Corollary 4 provides for any $p=\frac {2n}{n-1}$ an {\it explicit} construction of an `almost' $\Lambda_p$-set
which is not a $\Lambda_q$-set for $q\geq \frac {2n}{n-2}$.
No explicit  constructions of proper $\Lambda_p$-sets for $2<p<4$ seem to be known and their existence results from
probabilistic arguments (see [B2], [B4]).

In view of (3.10), Corollary 4 also provides explicit almost Euclidean subspaces of dimension $\sim N^{\frac 4p-1}$
in $\ell_N^p$, for $p$ of the form $\frac {2n}{n-1}, n\geq 4$ (while their maximal dimension is $\sim N^{\frac 2p}$
for $2<p<\infty$).
To be compared with the result from [G-L-R] on explicit almost Euclidean subspaces of $\ell_n^1$.

Returning to Proposition 3, we have more generally

\proclaim
{Proposition 3$'$}
Let $S$ be as in Proposition 3 and $T\in GL_n(\Bbb R), \Vert T\Vert>1$, an arbitrary invertible linear transformation.
Let $\Cal E=\Bbb Z^n \cap T(S)$.
Then, letting $p=\frac {2n}{n-1}$, we have the inequality
$$
\Big[\int_{\Bbb T^n}\Big|\sum_{z\in\Cal E} a_z e^{2\pi ix.z}\Big|^p dx\Big]^{\frac 1p} \ll \Vert T\Vert^\ve\Big(\sum_{x\in\Cal E}|a_z|^2\Big)^{\frac 12}.\tag
3.11
$$
\endproclaim

\noindent
{\bf Proof.}
Consider the set
$$
\Cal  E'= \{T^{-1} z; z\in\Cal E\}\subset S
$$
which elements are at least $\frac 1{\Vert T\Vert}$-separated.
Applying Proposition 2 with $D=1$ and $\rho=\frac 1{\Vert T\Vert}$, we obtain
$$
\lim_{R\to\infty} \Big|\nint_{B_R} \Big|\sum_{z\in\Cal E} a_z e^{2\pi i x'.T^{-1}z}\Big|^p dx'\Big] ^{\frac 1p}\ll\Vert T\Vert^\ve\Big(\sum_{z\in\Cal E}
|a_z|^2\Big)^{\frac 12}.\tag 3.12
$$

By change of variables $x= (T^{-1})^*x'$, it follows that
$$
\lim_{R\to\infty} \Big[\nint_{(T^{-1})^*(B_R)}\Big|\sum_{z\in \Cal E} a_z e^{2\pi i x.z}\Big|^p dx\Big]^{\frac 1p} \ll \Vert T\Vert^\ve\Big(\sum_{z\in\Cal E}
|a_z|^2 \Big) ^\frac 12\tag 3.13
$$
which, by periodicity, is equivalent to (3.11).

Take $S=\{(y, |y|^2); y\in\Bbb R^n, |y|<1\}$ the truncated paraboloid in $\Bbb R^{n+1}$ and let $T(x, t)=(Rx, R^2 t), R>1$.
From Proposition 3$'$, we immediately derive the following Strichatz' type inequality for the periodic Schr\"odinger group $e^{it\Delta}$.

\proclaim
{Corollary 5} Denote $\Delta$ the Laplacian on $\Bbb T^n$.
Then, for $p=\frac {2(n+1)}n$, we have the inequality
$$
\Vert e^{it\Delta} f\Vert_{L^p(\Bbb T^n \times \Bbb T)} \ll R^\ve \Vert f\Vert_{L^2(\Bbb T)}\tag 3.14
$$
assuming supp\,$\hat f \subset B(0, R)$.
\endproclaim

This bound should be compared with the following result established in [B3].

\proclaim
{Proposition 6} Let $f\in L^2(\Bbb T^n)$, $\Vert f\Vert_2 =1$ and such that supp\,$\hat f\subset B(0, R)$.
Then, for $\lambda>R^{\frac n4}$ and $q> \frac {2(n+2)}n$, the following inequality holds
$$
\mes [(x,t) \in\Bbb T^{n+1}; |e^{it\Delta} f|(x)>\lambda] < C_q R^{\frac n2 q-(n+2)}\lambda^{-q}.\tag 3.15
$$
\endproclaim

Combining Corollary 5, Proposition 6, we obtain the following improvement over Proposition 3.110 in [B3].

\proclaim
{Corollary 7} Let $n\geq 4$ (for $n<4$, better result may be obtained by arithmetical means, cf. [B3]).

Let $f$ be as in Proposition 6.
Then, for $q>\frac {2(n+3)}n$
$$
\Vert e^{it\Delta} f\Vert_{L^q(\Bbb T^{n+1})} < C_q R^{\frac n2 -\frac {n+2}q}\tag 3.16
$$
holds.
\endproclaim

Note that (3.16) is optimal.

\noindent
{\bf Proof.}

Denote $q_0=\frac {2(n+1)}n$ and $q_1$ some exponent $>\frac {2(n+2)}n$.
Let $F(x, t)=(e^{it\Delta} f)(x)$ and estimate for $q>q_1$
$$
\aligned
\int_{\Bbb T^{n+1}} |F|^q&\leq \int_{|F|>R^{\frac n4}} |F|^q+ R^{\frac n4(q-q_0)}\int|F|^{q_0}\\
&< C_{q_1} R^{\frac n2 q_1-(n+2)} \int^{R^{\frac n2}}_{R^{\frac n4}} \lambda^{q-1-q_1} d\lambda +C_\ve R^{\frac n4(q-q_0)+\ve}\\
&C_{q_1}\frac 1{q-q_1} R^{\frac n2 q-(n+2)}+C_\ve R^{\frac n4(q-q_0)+\ve} < C_q R^{\frac n2 q-(n+2)}
\endaligned
$$
for $q$ as above.

Corollary 5 admits a generalization that we discuss next.
Assume $\psi:\cup \to\Bbb R, U\subset\Bbb R^n$ a neighborhood of $0$, is a smooth function such that
$D^2\psi$ is positive (or negative) definite.
Then one has

\proclaim
{Proposition 8} Let $p=\frac {2(n+1)} n$ and $N\to\infty$.
Then for all $\ve>0$,
$$
\aligned
&\Big[\int_{[0, 1]^{n+1} \Big|\sum_{z\in\Bbb Z^n, \frac zN\in U} a_z e^{2\pi i}(x.z+N^2t\psi(\frac zN))}\Big|^pdxdt\Big]^{\frac 1p}\ll\\
&N^\ve\Big(\sum|a_z|^2\Big)^{1/2}.
\endaligned
\tag 3.17
$$
\endproclaim

Note that a coordinate change $x\mapsto x+Nt\nabla\psi(0)$ permits to assume $\psi(0)=\nabla\psi(0)=0$.
Let $S=[(x, \psi(x), x\in U]$ and
$$
\Cal E=\Big\{\Big(\frac zN, \psi\Big(\frac zN\Big)\Big); z\in\Bbb Z^n, \frac zN\in U\Big\}\subset S.
$$
Application of Proposition 2 with $\rho\sim\frac 1N$ implies that
$$
\aligned
&\Big[\int_{[0, 1]^{n+1}} \Big|\sum_{z\in\Bbb Z^n, \frac zN\in U} a_z e^{2\pi i(Nz.x+N^2\psi(\frac zN)t)}\Big|^p dxdt\Big]^{\frac 1p}\ll\\
&N^\ve\Big(\sum|a_z|^2\Big)^{1/2}
\endaligned
\tag 3.18
$$
and (3.17) follows by exploiting periodicity in $x$.
This proves Proposition 8.

Finally, observe that by taking $\psi(x) =\alpha_1x_1^2+\cdots+ \alpha_nx_n^2$ with $\alpha_1, \ldots,~\alpha_n~>~0$, Corollary 5 generalizes to a Strichartz inequality
for irrational tori, as considered in [B].
Applications to nonlinear Schr\"odinger type equations will not be discussed in this paper.


\Refs
\widestnumber\no{XXXXXX}

\ref\no{[B]}
\by J.~Bourgain
\paper On Strichartz inequalities and NLS on irrational tori
\jour Mathematical Aspects of Nonlinear Dispersive Equations, Annals of Math. 63 163 (2007), 1--20
\endref

\ref\no{[B1]}
\by J.~Bourgain
\paper Eigenfunction bounds for the Laplacian on the $n$-torus
\jour IMRN, 1993, no 3, 61--66
\endref

\ref\no{[B2]}\by J.~Bourgain
\paper
Bounded orthogonal systems and the $\Lambda(p)$-set problem
\jour Acta Math. 162 (1989), no 3, 227--245
\endref

\ref\no{[B3]}\by J.~Bourgain
\paper Fourier restriction phenomena for certain lattice subsets and applications to nonlinear evolution equations, I, Schr\"odinger equations
\jour GAFA 3 (1993), no 2, 107--156
\endref

\ref \no{[B4]}\by J.~Bourgain
\paper $\Lambda_p$-sets in analysis: results, problems and related aspects
\jour Handbook of the geometry of Banach Spaces, Vol I, 195--232
\endref

\ref\no{[B-C-T]} \by J.~Bennett, A.~Carbery, T.~Tao
\paper On the multilinear restriction and Kakeya conjectures
\jour Acta Math. 196 (2006), 202, 261--302
\endref

\ref\no{[B-G]} \by J.~Bourgain, L.~Guth
\paper Bounds on oscillatory integral operators obtained from multi-linear estimates
\jour preprint (to appear in GAFA)
\endref

\ref\no{[B-R-S]}\by J.~Bourgain, Z.Rudnick
\paper Restriction of toral eigenfunctions to hyper surfaces and nodal sets
\jour
(preprint)
\endref

\ref
\no{[G-L-R]}\by V.~Guruswami, J.~Lee, A.~Razborov
\paper Almost Euclidean subspaces of $\ell_1^N$ via expander codes
\jour Proc ACM-SIAM Symp. on Discrete Algorithms, 353--362, ACM (2008)
\endref





\endRefs
\enddocument